\documentclass[11pt]{article}
\usepackage{amssymb,amsthm,amsmath}
\usepackage[numbers,sort&compress]{natbib}
\usepackage{color}
\usepackage{graphicx}
\usepackage{tikz}
\def\eqref{[\ref]}

\newcommand\R{\mathbb R}

\newcommand{\bi}{\bibitem}

\begin{document}
\bigskip\bigskip\bigskip
\begin{center}
{\Large\bf %On point spectrum at critical coupling.
Alexander Gordon}\\
\bigskip\bigskip
%S. Jitomirskaya\footnote{Department of Mathematics, UC Irvine, USA}
\end{center}
\bigskip

Alexander (Sasha)  Gordon, died in Chicago on May 13, 2019, after a
 long
illness. Sasha was a brilliant mathematician, author of a number of
beautiful and original results in diverse fields of spectral theory
and other areas of analysis. His insights were behind some of the foundations of almost periodic
operators as well as the theme of genericity of singular continuous spectrum
developed by B. Simon and collaborators. Born in Kharkov, Ukraine, on
May 10, 1947, Sasha has had an unusual path in mathematics. At the time of his death, he was an associate
professor of mathematics at UNCC. His life story can be viewed as a
triumph of human and mathematical spirit over the circumstances.

Sasha is best known for what is called %the
 Gordon's lemma, a statement
about absence of decaying solutions of 1D equations close to periodic.
It has led to the first example of an almost periodic
operator with singular-continuous spectrum, to the understanding of
importance of arithmetics in this area, to ultimately sharp arithmetic spectral
transitions (e.g. \cite{ayz,jl,jk}) as well as powerful consequences in the
field of 1d quasicrystals (e.g. \cite{d}).

Before Sasha's work it was not even known whether a condition on the frequency,
imposed in various KAM arguments to prove point spectrum, is
necessary (it was sometime before the singular continuous revolution
of the 90s which allowed for an easy soft argument). His proof is
based on a beautiful in its simplicity Lemma on $SL(2,\R)$ matrices,
formulated and proved in his two-page ``Uspehi''\footnote{The top Russian journal at the time, publishing
  mainly announcements with the two-page limit.} paper \cite{gl}. Later
E.A. Gorin (1935-2018)\cite{gorin} found a sleeker proof that easily generalized
to a $GL(n,\R)$ version, using the
Cayley-Hamilton theorem, a result that remained unpublished. Gordon's paper appeared, in Russian, in 1975, and
was not known in the west until B. Simon's Moscow visit in 1981, when
S. Molchanov communicated to him, on the blackboard, Sasha's lemma with Gorin's proof.
Avron and Simon at the time were trying to implement Sarnak’s
suggestion that spectral properties of quasiperiodic operators might depend on arithmetic
properties of the frequencies, and this  lemma then led to their first
example of an operator
with singular-continuous spectrum: the super-critical almost Mathieu
operator with Liouville frequency \cite{as}.  This lemma was prominently featured in many
of Simon's articles and textbooks under the name ``Gordon's lemma'', which  Sasha
made several futile attempts to change to Gordon-Gorin's, see, for instance,
Footnote 2 in \cite{jst}. With the 1D result so natural and so simple, Sasha
went on a life-long quest to find a multi-dimensional version. There
have been attempts by various authors, not leading to anything
significant. Sasha finally succeeded recently, in a highly original joint work with
his friend from the undergraduate days A. Nemirovsky \cite{imrn}.
They  provided the first quantitative condition for absence of
point spectrum of multidimensional quasiperiodic operators, a result
he was still not fully satisfied with and thriving to make stronger till his last days.

Another fundamental work of comparable importance is Gordon's
theorem that 1D Schrodinger operators with interval spectrum have no
point spectrum for generic boundary values \cite{8,
9}. This indirectly influenced
B. Simon's discovery of  the Wonderland theorem \cite{won} and other results on
generic nature of singular continuous spectrum, a big theme in the
90s.

A certain prelude to that was a remarkable paper \cite{determ}
where Gordon constructed an explicit potential with Green's function
decay outside a set of energies of measure zero (leading to
localization for a.e. boundary value by spectral averaging). It motivated powerful results
by Kirsch, Molchanov, Pastur and others, as well as the work
described above because the
example, being so explicit, led
to a natural question whether the ``a.e.'' in the boundary values is
necessary - something that Sasha answered in a surprising and very
general way in \cite{8,9}.

Another fundamental contribution was the proof of measurability of
eigenelements of continuous operator families, first in a joint work
with Jitomirskaya, Last, and Simon \cite{gjls} and then in the work
with Kechris \cite{gk}. The issue is that for ergodic families, say,
the entire collection of eigenvalues is not a measurable object,
because it is invariant, yet not a.e. constant. Gordon's work showed,
for example,
that, in a very general setting, there exists a measurable
enumeration. This was fundamental for several further important
advances, particularly recent results on localization as a corollary
of dual reducibility for quasiperiodic operators, and first results on
arithmetic localization in the multi-dimensional setting.

Finally, we mention Sasha's proof  \cite{hp} of the Hartman-Putnam conjecture
that the lengths of spectral gaps of 1-D Schrodinger operators
with bounded potentials tend to zero at infinity, a long open problem. Hartman and Putnam
in their 1950 paper \cite{hp1}
proved this under an additional condition of uniform
continuity and asked whether it could be removed.

Sasha produced several gems also outside the spectral theory. The
highlights include
\begin{itemize}
\item A counterexample to the long-standing conjecture of A. Kolmogorov (1953) about the impossibility of mixed spectrum for an analytic flow with an integral invariant on a 2-D torus.\cite{kolm}

\item An effective sufficient condition for the cohomological equation
  to have no measurable solutions \cite{coh}
\end{itemize}

It should be noted that most of the beautiful results described above
were obtained when Sasha was doing math in his spare time.
His story is
in a way characteristic of time and circumstances. %, yet quite unique
Sasha got
his undergraduate degree from Moscow State University in 1970, with
high honors. As an undergraduate, he was an active participant of
the Banach Algebras
seminar led by E. Gorin and V. Lin, published two independent
papers in top Russian journals and wrote an excellent MS
thesis. Should he have finished in 1968 or earlier, such a
performance would have more than guaranteed him a recommendation to
continue as a PhD student, requiring only to pass several formal
examinations. Things had changed however in 1969, when 39 Jewish
students recommended for the PhD by their advisors, got Cs on the ``History of the
Communist Party'' exam and were not allowed in.
This marked the beginning of an era of significantly increased antisemitism at the School of Mathematics at
MSU and in Soviet math in general, an era that lasted till the end of
the 80s and of which Sasha was one of the victims. From 1970 on, Jews
were largely not even recommended to be admitted to the graduate admission exams, although
the problem soon almost disappeared because they mostly stopped being
admitted as undergrads already in 1968.\footnote{See,
      for instance, \cite{shif}
  for an account of the techniques used to achieve this in the society
  that proclaimed its ideals of equality and  internationalism, and
  for a fascinating story of an underground school, Jewish People's
  University, created in 1978 so that "Jewish children can learn
  math".}
% This marked the beginning of an era of renewed
% antisemitism at the math department of MSU and in Soviet mathematics in general, era
% that lasted till the end of the 80s and of which Sasha was one of the first
% victims.
% From 1970 on, Jews were largely not even recommended to be admitted to the
% exams, although the problem soon almost disappeared because they  mostly stopped being
% admitted as undergrads either.
Michael Brin, Svetlana Katok, Yakov Pesin are some other
names of mathematicians who finished their undergraduate studies in
1969-70 and were treated similarly. For an excellent account of related issues see
\cite{katok}.

For the next almost 25 years Sasha's work in math was done purely for
fun. He took a day  job but was able to stay in Moscow which made it
possible for him to attend research seminars at MSU. The MSU math department
was then at the beginning of the end of its “golden age”, and the
quality of  seminars was still outstanding (see again \cite{katok} for
more detail). Sasha regularly attended the one on Mathematical Physics led by
S. Molchanov, A. Ruzmaikin and D. Sokoloff. There, he was exposed to problems in spectral theory that led to his papers.

His advisor was S. Molchanov, whose questions proposed at the
seminar indeed stimulated much of Sasha's research in spectral theory,
but who views his own role rather as that of a friend and benefactor than
an advisor. Sasha's thesis featuring, among other things,
Hartman-Putnam and absence of point spectrum, was ready by about 1980,
but to get a PhD, it was necessary for the advisor to find a
university which would agree to schedule a defense for an outsider. It usually involved also finding a
fake advisor, thus providing no benefit to the academic record of the
true advisor. Most Jewish PhDs in the 70s-80s were obtained this
way. Yet it also took time and financial resources to travel and be away from
the day job, something that Sasha didn't have, making the task
especially complicated in his case. Molchanov's multiple attempts to
organize the defense at various universities failed. The university in
Sasha's native Kharkov that was successfully used for this purpose in
the 70s, say for the defense of M. Brin, by the 80s would no longer allow their
own Jews, with, for example, M. Lyubich having to defend in
Uzbekistan. Things changed around 1987, and in 1988  Molchanov
finally succeeded in organizing Sasha's defense at the Moscow Institute of  Electronic
Machine Building where  V. Maslov
(from Maslov's index)
was a  department chair and thought it would be easy,
but at the end, it still took Maslov threatening to quit his job to make it
happen. Finally, in 1993 Sasha got a job at a research laboratory that
would utilize his PhD and pay for him doing mathematics. At about the same
time, the Russian government
essentially stopped paying not just living wages, but any wages at all,
to
scientists.

As the Iron Curtain fell, Sasha got emboldened by the realization that his
earlier work, popularized by Molchanov and Simon, has found interest and
acclaim, and started trying to realize the dream of his youth to become
a professional mathematician. He had several short visiting positions
in Europe, and when invited by Molchanov to visit UNCC for a year in 1995-96, he
decided to try to make it in the US. Yet, despite the Gordon Lemma fame, for a 48
year old with a heavy accent and a short list of short publications full of 3-5 year gaps, this
seemed an almost impossible quest. There was not much he could do
about the age or the accent, but remarkably, in the next ten years he
managed to change the publication list around, producing a steady
record of publications in computational bio-statistics, other aspects
of applied math (true random number generators), and squeezing in a few
in his beloved spectral theory. His work in bio-statistics was
related to his job as a programmer at the Department of Computational
Biology, University of Rochester, 2001-2006. It has been published in
some of the leading journals in the field, including Annals of Applied
Statistics and SIAM J.Appl.Math, so we believe he has brought his
remarkable originality there as well. In 2006, Sasha's UNCC friends
managed to achieve what seemed to be impossible, and, aged 59, Sasha became a
tenure-track assistant professor there. He was happy and grateful
for an opportunity to live a normal academic life. He continued
working in both spectral theory and computational bio-statistics.

Sasha's papers were usually short and always struck with their
simplicity, clarity and mathematical beauty.
While his list of early publications is much shorter than it could
have been hadn't he had so many obstacles, it is fair to say that most
of his works from that period, those that found lots of resonance as well as those that went almost unnoticed outside Russia, are true mathematical masterpieces.

Outside mathematics, Sasha was a connoisseur of Russian literature,
especially poetry. He was active at UNCC's Russian literature club,
where he often led the readings.  He knew thousands of poems by heart
and also wrote beautiful poems himself.
% that he only shared with friends.

Sasha was battling a grave illness for the last six years, yet he
continued working and creating until the very end. He finished grading his students’ exams from
his hospital bed two days before his death. His
last papers, one each in computational bio-statistics and in the
spectral theory  are being published posthumously, the latter one in the
present issue. He will be greatly missed.\\

\vskip .2in

S. Jitomirskaya, S. A. Molchanov, B. Simon, B. R. Vainberg

\end{document}